\documentclass[12pt]{article}
\usepackage{amsmath, graphicx, amsfonts,amssymb}
\usepackage{color}

\def\be{\begin{equation}}
\def\ee{\end{equation}}
\def\bea{\begin{eqnarray}}
\def\eea{\end{eqnarray}}
\def\bt{\begin{theorem}}
\def\et{\end{theorem}}
\def\bl{\begin{lemma}}
\def\el{\end{lemma}}
\def\br{\begin{remark}}
\def\er{\end{remark}}
\def\bc{\begin{corollary}}
\def\ec{\end{corollary}}
\def\bd{\begin{definition}}
\def\ed{\end{definition}}

\def\b{\beta}

\def\k{\kappa}

\def\m{\mu}

\def\bbN{\mathbb{N}}

\def\bbR{\mathbb{R}}

\def\bn{B_2^}

\def\b1{B_{1}^}

\def\pt{\partial}

\def\ba{\begin{array}}
\def\ea{\end{array}}
\def\ben{\begin{enumerate}}
\def\een{\end{enumerate}}
\newtheorem{theorem}{Theorem}[section]
\newtheorem{lemma}{Lemma}[section]
\newtheorem{remark}{Remark}[section]

\newtheorem{corollary}{Corollary}[section]

\newtheorem{definition}{Definition}[section]
\begin{document}
\title{On the monotone properties of general affine surface areas under the Steiner
symmetrization }

\author{Deping Ye }
\date{}
\maketitle
\begin{abstract} In this paper, we prove that, if functions (concave) $\phi$ and
(convex) $\psi$ satisfy certain conditions, the $L_{\phi}$ affine
surface area is monotone increasing, while the $L_{\psi}$ affine
surface area is monotone decreasing under the Steiner
symmetrization. Consequently, we can prove related affine
isoperimetric inequalities, under certain conditions on $\phi$ and $\psi$, without assuming that the convex body involved has centroid (or the Santal\'{o}  point) at the origin.

\end{abstract}

\section{Introduction}

Affine invariants are very important tools for geometric analysis.
Many powerful affine invariants arise in the extension of the
Brunn-Minkowski theory -- started by Lutwak, for instance, the
$L_p$ affine surface areas introduced by Lutwak in the ground
breaking paper \cite{Lu1}. Notice that the study of the classical
affine surface area even went back to Blaschke \cite{Bl1} in 1923.
The $L_p$ affine surface areas have been proved to be key
ingredients in many applications, for instance, in the theory of
valuations (see e.g.  \cite{A1, A2, K, LudR, LR1}), approximation
of convex bodies by polytopes (see e.g., \cite{Gr2, LuSchW, SW5}),
and information theory (for convex bodies, see e.g.,
\cite{Jenkinson2012, Paouris2010, Werner2012a, Werner2012b}). A
beautiful result by Reisner, Sch\"{u}tt and Werner
\cite{ReisnerSchuttWerner2010} even implies that the Mahler volume
product (related to the famous unsolved Mahler conjecture) attains
the minimum only at those convex bodies with $L_p$ affine surface
areas equal to zero for all $p\in(0, \infty)$.

 Recently, much effort has been made to extend the $L_p$
Brunn-Minkowski theory to its next step: the
Orlicz-Brunn-Minkowski theory, which is of great demand, (see
e.g., \cite{Chen2011, HaberlLYZ,  Li2011, Ludwig2009,  LR1,
LYZ2010a, LYZ2010b, Deping2012, Zhu2012}). As mentioned in
\cite{LYZ2010b}, ``this need is not only motivated by compelling
geometric considerations (such as those presented in Ludwig and
Reitzner \cite{LR1}), but also by the desire to obtain Sobolev
bounds (see \cite{HaberlFranz2009}) of a far more general nature."
In particular, the $L_p$ affine surface areas were extended in
\cite{Ludwig2009, LR1}. As an example, here we give the definition
for the $L_{\phi}$ affine surface area. Let $Conc(0,\infty)$ be
the set of functions $\phi: (0,\infty)\rightarrow (0,\infty)$ such
that either $\phi$ is a nonzero constant function, or $\phi$ is
concave with $\lim _{t\rightarrow 0} \phi(t)=0$ and $\lim
_{t\rightarrow \infty} \phi(t)/t=0$ (in this case, set
$\phi(0)=0$). Note that the function $\phi$ is monotone
increasing. For $\phi \in Conc(0,\infty)$, {\it the $L_{\phi}$ affine
surface area of $K$} \cite{Ludwig2009, LR1} takes the following
form
\begin{equation*} as_{\phi}(K) = \int _{\partial K} \phi
\left(\frac{\k_K (y)}{\langle y, N_K(y)\rangle^{n+1}}\right)
\langle y, N_K(y)\rangle\,d\mu_K(y).\end{equation*} Here, $N_K(y)$
is an outer unit normal vector at $y$ to $\partial K$--the
boundary of $K$, $\k _K(y)$ is the Gaussian curvature at $y\in \pt
K$, and $\m _K$ denotes the usual surface area measure on $\pt K$.
The standard inner product on $\bbR^n$ is $\langle \cdot,
\cdot\rangle$ and it induces the Euclidian norm $\|\cdot\|$. The
$L_p$ affine surface area for $p\geq 0$ is corresponding to
$\phi(t)=t^{\frac{p}{n+p}}$. Here, for all $p\neq -n$, {\it the $L_p$
affine surface area of $K$} is defined as  (see e.g. \cite{Bl1, Lu1, SW5})
\begin{equation*}
as_{p}(K)=\int_{\partial K}\frac{\kappa_K(x)^{\frac{p}{n+p}}}
{\langle x,N_{ K}(x)\rangle ^{\frac{n(p-1)}{n+p}}} d\mu_{ K}(x),
\end{equation*}
and
\begin{equation*}
as_{\pm \infty}(K)=\int_{\partial K}\frac{\kappa _K (x)}{\langle
x,N_{K} (x)\rangle ^{n}} d\mu_{K}(x),
\end{equation*}
provided the above integrals exist. A fundamental result on the
$L_{\phi}$ affine surface area is the characterization theory of
upper--semicontinuous $SL(n)$ invariant valuation \cite{LR1}.
Namely, {\em every upper--semicontinuous, $SL(n)$ invariant
valuation vanishing on polytopes can be represented as an
$L_{\phi}$ affine surface area for some $\phi\in Conc(0,\infty)$.}

  (Affine) isoperimetric inequalities, such as, the
celebrated Blaschke-Santal\'{o} inequality, are of particular
importance to many problems in geometry. An affine isoperimetric
inequality relates two functionals associated with convex bodies
(or more general sets) where the ratio of the two functionals is
invariant under non-degenerate linear transformations. Affine
isoperimetric inequalities are arguably more useful than their
better known Euclidean relatives. For instance, the classical
affine isoperimetric inequality (see \cite{Sch}) gives an upper
bound for the classical affine surface area in terms of volume. This
inequality has very important applications in many other problems
(e.g. \cite{GaZ, Lu-O}).  The affine isoperimetric inequalities for $L_{\phi}$
affine surface areas were established in \cite{Ludwig2009}.
Namely, {\em among all convex bodies with fixed volume and with
centroids at the origin, the $L_{\phi}$ affine surface area attains
its maximum at ellipsoids.} For the case of $L_p$ affine
surface areas, such inequalities were already appear in \cite{Lu1,
WY2008}.  We refer readers to Section 2 and
\cite{Ludwig2009, LR1} for other general affine surface areas and their
properties.

  Note that these nice affine isoperimetric inequalities
were established by using Blaschke-Santal\'{o} inequality, which
requires the centroid (or the Santal\'{o} point) to be the
origin. Removing such a restriction is one of the main motivations
of this paper. Our main tool is the famous Steiner symmetrization, a powerful tool in
convex geometry. It is of particular importance in proving many
geometric inequalities (especially with either maximizer or
minimizer to be Euclidean balls and/or ellipsoids). To fulfill the
goals, one needs (1) to prove the monotone properties of objects
of interest under the Steiner symmetrization, and (2) to employ
the fact that each convex body will eventually approach to (in
Hausdorff metric) an origin-symmetric Euclidean ball by the
successive Steiner symmetrization.

  In this paper, we will study monotone properties of
general affine surface areas under the Steiner symmetrization, and
then provide a different proof for related affine isoperimetric
inequalities proved in the remarkable paper by Ludwig
\cite{Ludwig2009}.  In literature, 
some special cases have been studied. For instance, in \cite{Hug1995}, Hug proved that the 
classical affine surface area (with respect to $\phi(t)=t^{\frac{1}{n+1}}$) is monotone increasing 
under the Steiner symmetrization, and hence the classical affine isoperimetric inequality follows. When $K$ is smooth enough, the 
$L_{\pm \infty}$ affine surface area of $K$ (with respect to $\phi(t)=t$) is equal to the volume of $K^\circ$, the polar body of $K$. 
A well-known result in \cite{MeyerPajor1990} states that, if $K$ has centroid (or the Santal\'{o} point) at the origin, then the volume
of $K^\circ$ is 
smaller than or equal to the volume of the polar body of the Steiner symmetral of $K$ under any direction; 
hence the celebrated Blaschke-Santal\'{o} inequality follows. Readers are
referred to, e.g., \cite{Hug1995, MeyerPajor1990} (and
references therein) for related work for $p=1, \pm\infty$.

We summarize our main results on the $L_{\phi}$
affine surface areas as follows. For all $\xi \in S^{n-1}$, the
Steiner symmetral of $K$ with respect to $\xi$ is denoted by
$S_{\xi} (K)$. We also use $|K|$ to denote the volume of $K$.

\noindent {\bf Theorem A.} {\em Let $K\subset \bbR^n$
be a convex body with the origin in its interior, and $\phi\in
Conc(0, \infty)$. Assume that the function $F(t)=\phi(t^{n+1})$
for $t\in (0, \infty)$ is concave. 

\vskip 2mm \noindent (i). The
$L_{\phi}$ affine surface area is monotone increasing under the
Steiner symmetrization, i.e., $as_{\phi} (K)\leq
as_{\phi}(S_{\xi}(K))$ for all $\xi\in S^{n-1}$.

\vskip 2mm \noindent (ii). The affine isoperimetric inequality for the 
$L_{\phi}$ affine surface area holds true. That is,
$$as_{\phi}(K)\leq as_{\phi}(B_K),$$ where $B_K$ is the
origin-symmetric ball such that $|K|=|B_K|$. 

\noindent Moreover, if
$F(\cdot) $ is strictly concave and $K$ has positive Gaussian curvature almost everywhere, then equality holds if and
only if $K$ is an origin-symmetric ellipsoid.}

  This paper is organized as follows. Section 2 is for
notations and background on convex geometry, especially for
general affine surface areas. The main results will be proved in
Section 3. General references for convex geometry are \cite{Ga,
Sch}.

\section{Background and Notations}
  The setting will be the $n$-dimensional Euclidean space
$\bbR^n$. The standard orthonormal basis of $\bbR^n$ are $\{e_1,
\cdots, e_n\}$. When we write $\bbR^n=\bbR^{n-1}\times \bbR$, we
assume that $e_n$ is associated with the last factor, and hence,
\begin{equation}\label{projection:en} e_n^\perp=\{z\in \bbR^n: \langle z,
e_n\rangle=0\}=\{z\in \bbR^n: z_n=0\}.\end{equation}  
A set $L\subset e_n^{\perp}$ will be identified as a subset (and still denoted by $L$) of $\bbR^{n-1}$ (by deleting the last coordinate, which is $0$ always). A {\it convex
body} $K\subset\bbR^n$ is a compact convex subset of $\bbR^n$ with
nonempty interior. In this paper, we always assume that the origin
is in the interior of $K$.  We use $B^n_2$ and $S^{n-1}$ to denote
the unit Euclidean ball and sphere in $\bbR^n$ respectively.  The
{\it polar body} of $K$, denoted by $K^\circ$, is defined as
$K^{\circ}=\{y\in \bbR^n: \langle x,y \rangle \leq 1, \forall x\in
K\}.$ The {\it support function} $h_K$ of $K$ is defined as $h_K(u)=\max_{x\in K}\langle x, u\rangle$,
for all $u\in S^{n-1}$. The {\it Hausdorff distance}
between two convex bodies $K, L\subset \bbR^n$, denoted by $d_H(K,
L)$, is defined as
$$d_H(K, L)=\|h_K(u)-h_L(u)\|_{\infty}=\max_{u\in
S^{n-1}}|h_K(u)-h_L(u)|.$$ A convex body $K$ is said to have {\it curvature function} $f_K(u): S^{n-1}\rightarrow \bbR$ if 
$$ V(L, \underbrace{K, \cdots, K}_{n-1})=\frac{1}{n}\int_{S^{n-1}}h_L(u)f_K(u)\,d\sigma (u),$$  where 
$ V(L, {K, \cdots, K})$ is the mixed volume and $\sigma$ is the classical spherical measure on $S^{n-1}$. The {\it $L_p$ curvature function} for $K$ is defined as $f_p(K, u)=h_K(u)^{1-p}f_K(u)$ (see \cite{Lu1}).

For a linear map or a matrix $T$, we use $\det (T)$ for
the determinant of $T$. For a (smooth enough) function $f:
\bbR^{n-1}\rightarrow \bbR$, $\nabla f(x)$ denotes its gradient
function, and  $\langle f(x)\rangle=f(x)-\langle
x, \nabla f(x)\rangle$. Note that $\langle f(x)\rangle$ is linear;
namely,  for any two (smooth enough) functions $f, g:
\bbR^{n-1}\rightarrow \bbR$ and for all $a, b\in \bbR$, one has
$\langle af(x)+bg(x)\rangle=a \langle f(x)\rangle+b\langle
g(x)\rangle$. In particular, we will often use the following
special case $$\big\langle
\frac{f(x)+g(x)}{2}\big\rangle=\frac{ \langle
f(x)\rangle+\langle g(x)\rangle}{2}.$$

In \cite{Ludwig2009}, Ludwig also introduced the
$L_{\psi}$ affine surface area. Hereafter, $Conv(0,\infty)$
denotes the set of functions $\psi: (0,\infty)\rightarrow
(0,\infty)$ such that either $\psi$ is a nonzero constant
function, or $\psi$ is convex with $\lim_{t\rightarrow 0}
\psi(t)=\infty$ and $ \lim_{t\rightarrow \infty} \psi(t)=0$ (in
this case, we set $\psi(0)=\infty$). For $\psi \in
Conv(0,\infty)$, {\it the $L_{\psi}$ affine surface area of $K$} can be
formulated as
\begin{equation*}as_{\psi}(K) =
\int _{\partial K} \psi \left(\frac{\k_K (y)}{\langle y,
N_K(y)\rangle^{n+1}}\right) \langle y, N_K(y)\rangle\,d\mu_K(y)
.\end{equation*} In particular, the $L_p$ affine surface area for
$-n<p<0$ is corresponding to $\psi(t)=t^{\frac{p}{n+p}}$. Moreover,
for all $\psi\in Conv(0, \infty)$, the $L_{\psi}$ affine surface
area is a lower--semicontinuous, $SL(n)$-invariant valuation.
Affine isoperimetric inequalities for $L_{\psi}$ affine surface
areas were also established in \cite{Ludwig2009}. Namely, {\em
among all convex bodies with fixed volume and with centroids at the
origin, the $L_{\psi}$ affine surface area attains its minimum
 at ellipsoids.}

 The $L_{\phi}$ affine surface area for $\phi\in Conc(0,
\infty)$ is proved to be upper--semicontinuous \cite{LR1}. That is,
for any sequence of convex bodies $K_i$ converging to $K$ in the
Hausdorff metric $d_H(\cdot, \cdot)$, one has $$\lim
\sup_{j\rightarrow \infty} as_{\phi} (K_j)\leq as_{\phi} (K). $$
The $L_{\psi}$ affine surface area for $\psi\in Conv(0, \infty)$
is showed to be lower--semicontinuous \cite{Ludwig2009}; namely for
any sequence of convex bodies $K_i$ converging to $K$ in the
Hausdorff metric $d_H(\cdot, \cdot)$, one has $$\lim
\inf_{j\rightarrow \infty} as_{\psi} (K_j)\geq as_{\psi} (K). $$
The semicontinuous properties will be crucial in proving related
affine isoperimetric inequalities in Section 3.

For a convex body $K\subset \bbR^n$ and on the
direction $e_n$,  let $K_{e_n}\subset e_n^\perp$ be the orthogonal projection of $K$ onto $e_n^\perp$ and $K_0\subset e_n^{\perp}$ be the relative interior of $K_{e_n}$. We denote $f, g: K_{e_n} \rightarrow \bbR$ the
{\it overgraph} and {\it undergraph} functions of $K$ with respect to $e_n$,
i.e.,
\begin{eqnarray*}  
f(x)=max\{t\in \bbR: (x,t)\in K\},  \ 
g(x)=-min\{t\in \bbR: (x,t)\in K\},\    x\in K_{e_n}.
\end{eqnarray*} 
Denote $K^{+}:=graph(f(x))$ and $K^{-}:=graph(-g(x))$. Similarly, one can define overgraph and undergraph functions for $K$ with respect to any direction $u\in S^{n-1}$.

\bl \label{Equivalent:definition:affine:surface:area} Let
$K\subset \bbR^n$ be a convex body with the origin in its
interior.  For all $\phi\in Conc(0, \infty)$, one has
$$as_{\phi}(K)=\int _{K_0}\left\{\phi\left(\frac{|\mathrm{det}
(d^2f(x))|}{\langle f(x)\rangle^{n+1}}\right) \langle f(x)\rangle
+\phi\left(\frac{|\mathrm{det} (d^2g(x))|}{\langle
g(x)\rangle^{n+1}}\right) \langle g(x)\rangle \right\}\,dx.$$
 \el

 \noindent {\bf Proof.} For almost all $x\in K_0$ (with
respect to the $(n-1)$-dimensional Lebesgue measure on $K_0$), the
Gaussian curvature of the point $y=(x, f(x))\in
\partial K$ can be formulated as \cite{Lei1986}
\begin{equation}
\kappa _K(y)=\frac{|\mathrm{det}(d^2f(x))|}{(\sqrt{1+\|\nabla
f(x)\|^2})^{n+1}}. \label{Gauss:curvature:represation}
\end{equation}  On the other hand, the outer
unit normal vector to $\partial K$ at the point $y$ is             
\begin{equation}\label{Outer:normal:vector}
N_K(y)=\frac{(-\nabla f(x), 1)}{\sqrt{1+\|\nabla f(x)\|^2}}.
\end{equation}
Hence, we have
\begin{equation}
\label{Support:function:representation} \langle y,
N_K(y)\rangle=\frac{f(x)-\langle x, \nabla
f(x)\rangle}{\sqrt{1+\|\nabla f(x)\|^2}}=\frac{\langle
f(x)\rangle}{\sqrt{1+\|\nabla f(x)\|^2}}.
\end{equation}
For almost all $x\in K_0$, by formulas
(\ref{Gauss:curvature:represation}), (\ref{Outer:normal:vector}),
and (\ref{Support:function:representation}), one has, at the point
$y=(x, f(x))\in \partial K$,
\begin{eqnarray}
\frac{\kappa_K(y)} {\langle y,N_{K}(y)\rangle
^{n+1}}&=&\frac{|\mathrm{det} (d^2f(x))|}{\langle
f(x)\rangle^{n+1}}. \label{representation:affine:invariant}
\end{eqnarray}
The surface area can be rewritten as
\begin{equation*}
\,d\mu_K(y)=\sqrt{1+\|\nabla f(x)\|^2}\,dx. \end{equation*}  
Combining with the Federer's area formula (see \cite{Federer1969}) and $f(x)$ being  locally Lipschitz,  one has
\begin{equation*}
\int_{K^+}\phi \left(\frac{\k_K (y)}{\langle y, N_K(y)\rangle^{n+1}}\right)
\langle y, N_K(y)\rangle\,d\mu_K(y)=\int_{K_0}\phi\left(\frac{|\mathrm{det}
(d^2f(x))|}{\langle f(x)\rangle^{n+1}}\right) \langle
f(x)\rangle\,dx.
\end{equation*}

\noindent Similarly, $g(x)$ being locally Lipschitz and the  Federer's area formula imply that 
 \begin{equation*} \int_{K^-}\phi
\left(\frac{\k_K (z)}{\langle z, N_K(z)\rangle^{n+1}}\right)
\langle z, N_K(z)\rangle\,d\mu_K(z)=\int_{K_0}\phi\left(\frac{|\mathrm{det}
(d^2g(x))|}{\langle g(x)\rangle^{n+1}}\right) \langle
g(x)\rangle\,dx.
\end{equation*}
Finally, let $K'=\partial K\cap (\mathrm{relbd} (K_0), \bbR)$ where $$(\mathrm{relbd} (K_0) ,\bbR)=\{(x,t): x\in \mathrm{relbd} (K_0), \ t\in \bbR\}.$$ Then, the boundary of $K$ can be decomposed as $K^+\cup K^-\cup  K'$. From generalized cylindrical coordinates, one gets, for $\phi\in Conc(0,\infty)$,
\begin{eqnarray*}
\int _{ K'} \phi
\left(\frac{\k_K (y)}{\langle y, N_K(y)\rangle^{n+1}}\right)
\langle y, N_K(y)\rangle\,d\mu_K(y)=0. \end{eqnarray*} Therefore, 
\begin{eqnarray*}
as_{\phi}(K)&=&\int _{ K^+\cup K^-} \phi
\left(\frac{\k_K (y)}{\langle y, N_K(y)\rangle^{n+1}}\right)
\langle y, N_K(y)\rangle\,d\mu_K(y)\\ &=&\int _{K_0}\left\{\phi\left(\frac{|\mathrm{det}
(d^2f(x))|}{\langle f(x)\rangle^{n+1}}\right) \langle f(x)\rangle
+\phi\left(\frac{|\mathrm{det} (d^2g(x))|}{\langle
g(x)\rangle^{n+1}}\right) \langle g(x)\rangle \right\}\,dx. \end{eqnarray*}

To have similar results for the $L_{\psi}$ affine surface area, one needs to assume that $K$ has positive Gaussian curvature almost everywhere (with respect to the measure $\mu_K$).  Along the same line, one can prove the following lemma.

\bl \label{Equivalent:definition:affine:surface:area:psi} Let
$K\subset \bbR^n$ be a convex body with the origin in its
interior and having positive Gaussian curvature almost everywhere.  For all $\psi\in Conv(0, \infty)$, one has
$$as_{\psi}(K)=\int _{K_0}\left\{\psi\left(\frac{|\mathrm{det}
(d^2f(x))|}{\langle f(x)\rangle^{n+1}}\right) \langle f(x)\rangle
+\psi\left(\frac{|\mathrm{det} (d^2g(x))|}{\langle
g(x)\rangle^{n+1}}\right) \langle g(x)\rangle \right\}\,dx.$$\el

\section{Main Results}
  The Steiner symmetral of $K$ with respect to
$e_n$, denoted by $S_{e_n}(K)$, is defined as
\begin{equation}\label{symmetrization:en}
S_{e_n}(K)=\left\{(x, t): -\frac{1}{2}(f(x)+g(x))\leq t\leq
\frac{1}{2}(f(x)+g(x)), x\in K_{e_n}\right\},\end{equation} where $K_{e_n}$ is the orthogonal projection of $K$ onto $e_n^\perp$. The Steiner symmetral of $K$ with
respect to any direction $\xi\in S^{n-1}$ is denoted by $S_{\xi} (K)$ and can be formulated
similar to formula (\ref{symmetrization:en}).  Clearly, the
Steiner symmetrization does not change the volume; namely, $|S_{\xi}(K)|=|K|$ for all
$\xi\in S^{n-1}$. In the later proof,
we focus on the direction $e_n=(0, \cdots, 0, 1)$ only.

 The following well-known lemma is of particular
importance in applications (see e.g. \cite{Ga}). \bl
\label{limit:symmetrization} Let $K$ be a convex body in $\bbR^n$.
There is a sequence of directions $\xi _i\in S^{n-1}$, $i\in
\bbN$, such that the successive Steiner symmetral of $K$
$$K_m=S_{\xi_m}(S_{\xi_{m-1}} (\cdots (S_{\xi_1}(K))\cdots ))$$
converges to an origin-symmetric Euclidean ball in Hausdorff
distance $d_H(\cdot, \cdot)$. \el

  The following lemma is an easy consequence of Theorem G
in \cite{Roberts1973} (see p. 205).

\bl\label{curvature:symmetrization:1} Let $A$ and $B$ be two
$(n-1)\times (n-1)$ symmetric, positive semi-definite matrices.
Then
$$2\left[\mathrm{det} \left(\frac{A+B}{2}\right) \right]^{\frac{1}{n+1}}\geq \big(\mathrm{det}(A) \big)^{\frac{1}{n+1}}
+\big(\mathrm{det}(B) \big)^{\frac{1}{n+1}}.$$  If, in addition,
$B$ (or $A$) is positive definite, equality holds if and only if
$A=B$. \el

In order to settle the equality conditions for general affine isoperimetric inequalities, we need the following lemma 
(see \cite{Hug1995}), which follows from Brunn's classical characterization of ellipsoids. For a convex body $K$ with the origin in its interior, and $\xi \in S^{n-1}$, we define by $M(K, \xi)$ the set of the midpoints of all line segments $K\cap L$ where $L$ varies over all lines with direction $\xi$ that meet the interior of $K$. 

\bl \label{Charactristic:Ellipsoid:Hug} Let $K$ be a convex body with the origin in its interior. Let $S^*$ be a dense subset of $S^{n-1}$. Then 
$K$ is an ellipsoid if and only if for each $\xi \in S^*$, the set $M(K, \xi)$ is contained in a hyperplane. \el 

The following lemma was proved in \cite{Hug1995}. 
\bl \label {Lemma 4.3} Let $U\subset \bbR^{n-1}$, $0\in U$, be open and convex. Let $f: U\rightarrow \bbR$ be locally Lipschitz and differentiable at $0$. If, $\langle x, \nabla f(x)\rangle=f(x)$ for almost all $x\in U$ such that $f$ is differentiable at $x$, then $f(x)=\langle v, x\rangle$ for all $x\in U$ and some suitable $v\in \bbR^{n-1}$. \el

\subsection{$L_{\phi}$ affine surface areas are increasing under  the
Steiner symmetrization}
\bt\label{symmetrization:L:phi} Let $K\subset \bbR^n$ be a convex
body with the origin in its interior and $\phi\in Conc(0,
\infty)$. Suppose that the function $F(t)=\phi(t^{n+1})$ for $t\in
(0, \infty)$ is concave. One has, for all
$\xi\in S^{n-1}$,
\begin{equation*}as_{\phi} (K)\leq as_{\phi}(S_{\xi}(K)).
\end{equation*}  \et 
{\bf Remark.}  Note that the function $\phi \in Conc(0,\infty)$ is monotone
increasing, hence $F$ is also an increasing function.   In fact, in view of Lemma \ref{curvature:symmetrization:1}, the condition that $F$ is concave and monotone increasing is more natural for proving Theorem \ref{symmetrization:L:phi}. If (non-constant) function $F$ is monotone increasing and concave with $F(0)=0$, then $\phi \in Conc(0,\infty)$. (It is easily checked that $F$ constant implies $\phi$ constant).  To this end, for all $0<t<s$, $\phi(t)=F(t^{\frac{1}{n+1}})\leq F(s^{\frac{1}{n+1}})=\phi(s)$ and hence $\phi$ is monotone increasing. Note that the function $t^{\frac{1}{n+1}}$ is concave, and by $F$ being increasing,  for all $\lambda\in [0,1]$ and $0<t<s$, 
\begin{eqnarray*} \phi(\lambda t+(1-\lambda) s)&=&F([\lambda t+(1-\lambda) s]^{\frac{1}{n+1}})\geq F(\lambda t^{\frac{1}{n+1}}+(1-\lambda) s^{\frac{1}{n+1}})\\ &\geq & \lambda F(t^{\frac{1}{n+1}})+(1-\lambda) F(s^{\frac{1}{n+1}})=\lambda \phi(t)+(1-\lambda)\phi(s).\end{eqnarray*}  
The concavity of $F$ implies that, for all (given) $t>1$ and $\lambda =1/t\in [0,1]$, $$F(1)=F(\lambda t+(1-\lambda)\cdot 0)\geq \lambda F(t)+(1-\lambda)F(0)=\frac{F(t)}{t}.$$ Hence, $F(t)\leq F(1)t$ for all $t>1$, which further implies that $$0\leq \lim_{t\rightarrow \infty} \frac{\phi(t)}{t}
=\lim_{t\rightarrow \infty} \frac{F(t^{\frac{1}{n+1}})}{t}\leq \lim_{t\rightarrow \infty} F(1)t^{\frac{-n}{n+1}}=0.$$
That is, $\lim_{t\rightarrow \infty}\phi(t)/t=0$. In words, we have conclude that $\phi\in Conc(0,\infty)$. However,  $\phi\in Conc(0,\infty)$ {\it does not} in general imply $F$ being concave and monotone increasing. For instance, for homogeneous function $\phi(t)=t^a$, to have $F$ concave, one needs $0\leq a\leq \frac{1}{n+1}$.

 \noindent {\bf Proof of Theorem \ref{symmetrization:L:phi}.} Let $K$ be a convex body with
the origin in its interior. Recall that, $as_{\phi}(K)$ is affine
invariant, i.e., for all
 linear maps with $|\det (T)|=1$ and for all $\phi\in Conc(0,\infty)$,
 one has
\begin{equation*} as_{\phi}(TK)=as_{\phi}(K).
\end{equation*}
Hence, without loss of generality, we only prove Theorem
\ref{symmetrization:L:phi} on the direction $$\xi =e_n=(0, \cdots,
0 ,1).$$  By Lemma
\ref{Equivalent:definition:affine:surface:area}, the $L_{\phi}$
affine surface area of $K$ can be rewritten as
$$as_{\phi}(K)=\int _{K_0}\left\{\phi\left(\frac{|\mathrm{det}
(d^2f(x))|}{\langle f(x)\rangle^{n+1}}\right)\langle f(x)\rangle
+\phi\left(\frac{|\mathrm{det} (d^2g(x))|}{\langle
g(x)\rangle^{n+1}}\right)\langle g(x)\rangle \right\}\,dx.$$
  Let $h(x)=[f(x)+g(x)]/2$.  By Lemma \ref{curvature:symmetrization:1}, one
gets, for almost all $x\in K_0$,
\begin{eqnarray} 2 \left|\mathrm{det}
\big(d^2h(x)\big)\right|^{\frac{1}{n+1}}&=& 2 \left|\mathrm{det}
\bigg[\frac{d^2f(x)+d^2g(x)}{2}\bigg]\right|^{\frac{1}{n+1}}\nonumber\\
&\geq& \left|\mathrm{det} (d^2 f(x))\right|^{\frac{1}{n+1}}
+\left|\mathrm{det} (d^2
g(x))\right|^{\frac{1}{n+1}}.\label{determinant:inequality:01}
\end{eqnarray}
Recall that $F(t)=\phi(t^{n+1})$ and $$\langle h(x)\rangle
=\frac{\langle f(x)\rangle +\langle g(x)\rangle }{2},$$ for almost all
$x\in K_0$. Hence,
\begin{eqnarray} \phi\left(\frac{|\mathrm{det}
(d^2h(x))|}{\langle h(x)\rangle^{n+1}}\right)=
F\left(\frac{|\mathrm{det} (d^2h(x))|^{\frac{1}{n+1}}}{\langle
h(x)\rangle}\right)\nonumber = F\left(\frac{2|\mathrm{det}
(d^2h(x))|^{\frac{1}{n+1}}}{\langle f(x)\rangle+\langle
g(x)\rangle}\right).\nonumber \end{eqnarray} Note that $\phi\in
Conc(0, \infty)$ is an increasing function, so is $F(t)$ on $(0,
\infty)$. By inequality (\ref{determinant:inequality:01}), one has
\begin{eqnarray}
\phi\left(\frac{|\mathrm{det} (d^2h(x))|}{\langle
h(x)\rangle^{n+1}}\right)&\geq& F\left(\!\frac{\left|\mathrm{det}
(d^2 f(x))\right|^{\frac{1}{n+1}} \!+\!\left|\mathrm{det} (d^2
g(x))\right|^{\frac{1}{n+1}}}{\langle f(x)\rangle+\langle
g(x)\rangle}\!\right)\nonumber\\
&\geq& F\left(\frac{|\mathrm{det} (d^2
f(x))|^{\frac{1}{n+1}}}{\langle f(x)\rangle} \right)\frac{\langle
f(x)\rangle}{\langle f(x)\rangle+\langle g(x)\rangle}\nonumber\\
&& \ + \ \ F\left(\!\frac{\left|\mathrm{det} (d^2
g(x))\right|^{\frac{1}{n+1}}}{\langle
g(x)\rangle}\!\!\right)\frac{\langle g(x)\rangle}{\langle
f(x)\rangle+\langle g(x)\rangle}, \label{concavity:inequality:1}
\end{eqnarray} where the last inequality in (\ref{concavity:inequality:1})
follows from the
concavity of the function $F(t)$ on $(0, \infty)$. Therefore, by
(\ref{symmetrization:en}) and Lemma
\ref{Equivalent:definition:affine:surface:area}, we have for all
$\phi\in Conc(0, \infty)$,
\begin{eqnarray*}
as_{\phi}(S_{e_n}(K))&=&2 \int
_{K_0}\left\{\phi\left(\frac{|\mathrm{det} (d^2h(x))|}{\langle
h(x)\rangle^{n+1}}\right)
\langle h(x)\rangle\right\}\,dx \nonumber\\
&\geq& \int _{K_0}\left\{\phi\left(\frac{|\mathrm{det}
(d^2f(x))|}{\langle f(x)\rangle^{n+1}}\right) \langle
f(x)\rangle+\phi\left(\frac{|\mathrm{det} (d^2g(x))|}{\langle
g(x)\rangle^{n+1}}\right) \langle g(x)\rangle\right\}\,dx
\nonumber
\\&=&as_{\phi}(K). \end{eqnarray*}

 Let $K$ be a convex body having curvature function $f_K: S^{n-1}\rightarrow \bbR$. For $\phi\in Conc(0, \infty)$, the
$L_{\phi}^*$ affine surface area of $K$, denoted by $as_{\phi}^*(K)$, can be
formulated as  \cite{Ludwig2009} 
$$as_{\phi}^*(K)=\int_{S^{n-1}}\phi(f_{-n}(K,u))\,d\nu_K(u),$$
where $f_{-n}(K,u)=h_K(u)^{n+1}f_K(u)$ is the $L_p$ curvature
function of $K$ (see \cite{Lu1}) for $p=-n$, while $\,d\nu_K(u)=\,d\sigma
(u)/h_K(u)^n$ with $\,d\sigma(u)$ the classical spherical measure
over the sphere $S^{n-1}$. It was proved in \cite{Ludwig2009} that
\begin{equation}\label{duality:concave}
as_{\phi}^*(K)=as_{\phi}(K^\circ),\end{equation} for any convex
body $K$ having curvature function and with the origin in its interior.
Combining with Theorem \ref{symmetrization:L:phi}, one immediately
has the following result.

\bc Let $K\subset \bbR^n$ be a convex body with the origin in its
interior and having curvature function. Let $\phi\in Conc(0, \infty)$. Assume that the function
$F(t)=\phi(t^{n+1})$ for $t\in (0, \infty)$ is concave. Then, the
$L_{\phi}^*$ affine surface area is monotone increasing under the
Steiner symmetrization in the following sense:
$$as_{\phi}^* (K)\leq as_{\phi}^*([S_{\xi}(K^\circ)]^\circ),$$ for
all $\xi\in S^{n-1}$, such that, $[S_{\xi}(K^\circ)]^\circ$ has curvature function. \ec


\bt \label{L:phi:affine isoperimetric inequality} Let $K\subset
\bbR^n$ be a convex body with the origin in its interior, and
$B_K$ be the origin-symmetric Euclidean ball with $|K|=|B_K|$.
Then, for all $\phi \in Conc(0, \infty)$ such that the function
$F(t)=\phi(t^{n+1})$ for $t\in (0, \infty)$ is concave, one has
\begin{equation*}as_{\phi}(K)\leq
as_{\phi}(B_{K}).\end{equation*} If in addition
$F(\cdot)$ is strictly concave and $K$ has positive Gaussian curvature almost everywhere, equality holds if and only if $K$
is an origin-symmetric ellipsoid. \et

\noindent {\bf Remark.} Clearly, if $as_{\phi}(K)=as_{\phi}(B_K)$, then $K$ {\it cannot} be  a convex body with Gaussian curvature equal to $0$ almost everywhere (with respect to the measure $\mu_K$) on the boundary of $K$, as otherwise $as_{\phi}(K)=0$.

 \noindent {\bf Proof.} Let $K\subset \bbR^n$ be a
convex body with the origin in its interior. Suppose that $\phi\in
Conc(0, \infty)$. By Lemma \ref{limit:symmetrization}, one can find  a sequence of directions   $\{u_i\}_{i=1}^{\infty}\subset \Omega$ such that  $K_i$ converges to
$B_K$ in the Hausdorff distance.  Here $K_i$ is defined as follows:  
$$K_1=S_{u_1}(K); \ \ \  K_{i+1}=S_{u_{i+1}}(K_i),\ \ \  \forall i=1, 2, \cdots $$ 
 Theorem \ref{symmetrization:L:phi} implies
that
$$as_{\phi}(K)\leq as_{\phi}(K_1)\leq \cdots \leq as_{\phi}(K_j), \ \ \forall
j\in \bbN. $$ Combining with the upper--semicontinuity of
$as_{\phi}(\cdot)$, one has,
\begin{eqnarray} as_{\phi}(K)&\leq& \lim
\sup _{j\rightarrow \infty} as_{\phi}(K_j)\leq as_{\phi}(\lim
_{j\rightarrow \infty}
K_j)=as_{\phi}(B_K).\label{Proof:Isoperimetric:phi}
\end{eqnarray}

 Now let us assume that $F$ is strictly concave and $as_{\phi}(K)=as_{\phi}(B_K)$.  Let $K$ be a convex body with positive  Gaussian curvature almost everywhere.
We now claim that the set $M(K, e_n)$ is contained in a hyperplane. In this case, we assume that $e_n$ is a direction such that  both the ovegraph and undergraph functions $f, g$ are differentiable  at $0$. Note that $as_{\phi}(K)=as_{\phi}(S_{e_n}(K))$ requires equalities for (\ref{concavity:inequality:1}). Combining with  Lemma \ref{curvature:symmetrization:1} and the strict concavity of $F$, one has, for almost every $x\in  K_0$, \begin{equation}\label{representation:affine:invariant-0}\left|\mathrm{det}
(d^2 f(x))\right|=\left|\mathrm{det} (d^2
g(x))\right|;  \ \ \  \frac{|\mathrm{det} (d^2
f(x))|^{\frac{1}{n+1}}}{\langle f(x)\rangle} =\frac{\left|\mathrm{det} (d^2
g(x))\right|^{\frac{1}{n+1}}}{\langle
g(x)\rangle}>0.\end{equation}  Hence, for almost all $x\in K_0$, 
$$\langle f(x)\rangle=f(x)-\langle x, \nabla f(x)\rangle=g(x)-\langle x, \nabla g(x)\rangle=\langle g(x)\rangle.$$ That is, $f(x)-g(x)=\langle x, \nabla (f(x)-g(x))\rangle$ for almost all $x\in K_0$. Note that $f-g$ is locally Lipschitz. From Lemma \ref{Lemma 4.3}, one obtains that $f(x)-g(x)$ is linear, and hence $M(K, e_n)$ is contained in a hyperplane. 

 Let $\Omega$ be the dense subset of $S^{n-1}$ such that the corresponding  overgraph and undergraph functions at the direction $u\in \Omega$ are both differentiable at $0$. Here, we have used the fact that $\sigma(S^{n-1}\setminus \Omega)=0$, because $u\in \Omega$ if and only if the radial function $\rho_K$ of $K$ is differentiable at $\pm u$.  For any $u\in \Omega$, there is a rotation $T$ such that $T(u)$ is parallel to $e_n$. The above claim then implies that $M(TK, T(u))$ (and hence $M(K, u)$) is contained in a hyperplane. By Lemma  \ref{Charactristic:Ellipsoid:Hug}, $K$ is an ellipsoid.  Moreover, $K$ has to be an origin-symmetric ellipsoid. To this end, without loss of generality, we assume that $K$ is a ball with center $y_c\neq 0$ (this follows from affine invariance of the $L_{\phi}$ surface area). By formulas (\ref{representation:affine:invariant}) and (\ref{representation:affine:invariant-0}), one gets 
 $$\frac{[\kappa_K(y)]^{\frac{1}{n+1}}} {\langle y,N_{K}(y)\rangle
}=\frac{|\mathrm{det} (d^2
f(x))|^{\frac{1}{n+1}}}{\langle f(x)\rangle} =\frac{\left|\mathrm{det} (d^2
g(x))\right|^{\frac{1}{n+1}}}{\langle
g(x)\rangle}=\frac{[\kappa_K(z)]^{\frac{1}{n+1}}} {\langle z,N_{K}(z)\rangle},\ \ \  \forall x\in K_{e_n},$$ with $y=(x, f(x))\in \partial K$ and $z=(x, -g(x))\in \partial K$. Note that the curvature of $K$ is a constant, and hence $\kappa_K(y)=\kappa_K(z)$. This implies  $\langle y,N_{K}(y)\rangle=\langle z,N_{K}(z)\rangle$, a contradiction with $K$ being a ball with center $y_c\neq 0$.

 Combining with formula (\ref{duality:concave}), one
immediately has the following isoperimetric inequality for the
$L_{\phi}^*$ affine surface area. \bc Let $K\subset \bbR^n$ be a
convex body with the origin in its interior and having curvature function. Let $B_K$ be the
origin-symmetric Euclidean ball with $|K|=|B_K|$. Then, for all
$\phi \in Conc(0, \infty)$ such that the function
$F(t)=\phi(t^{n+1})$ for $t\in (0, \infty)$ is concave, one has
\begin{equation*}as_{\phi}^*(K)\leq
as_{\phi}^*([B_{K^\circ}]^\circ).\end{equation*} If in addition $F(\cdot)$ is strictly concave and $K^\circ$ has positive Gaussian curvature almost everywhere (with respect to $\mu_{K^\circ}$), equality holds if
and only if $K$ is an origin-symmetric ellipsoid.\ec

 As the $L_p$ affine surface areas for $p>0$ are special
cases of $L_{\phi}$ affine surface areas, with
$\phi(t)=t^{\frac{p}{n+p}}$, we get the following result.
 \bc \label{Lp:p:positive:symmetrization} Let $K$ be a convex body with the origin in its interior, and
 let $p\in (0, 1)$.\\
 \noindent (i) The $L_p$ affine surface area for $p\in (0,1)$ is monotone increasing under the Steiner symmetrization.
 That is,  for any $\xi \in S^{n-1}$, one has
\begin{eqnarray*}as_p(K)\leq as_p(S_{\xi}(K)).
 \end{eqnarray*} (ii) The $L_p$
affine surface areas attain their maximum at the ellipsoid, among all convex bodies with fixed volume. More
precisely, \begin{eqnarray*} \frac{as_p(K)}{as_p(B^n_2)}\leq
\left(\frac{|K|}{|\bn n|}\right)^{\frac{n-p}{n+p}}. 
\end{eqnarray*} For convex bodies with positive Gaussian almost everywhere, equality holds if and only if $K$ is an origin-symmetric ellipsoid. \ec

 \noindent {\bf Remark.} Notice that if $p=0$, the $L_p$
affine surface area is equal to the volume and hence will not
change under the Steiner symmertrization. The case $p=1$
corresponds to the classical affine surface area, and this has
been proved, for instance, in \cite{Hug1995}. The $L_p$
affine isoperimetric inequalities for $p>1$ were first established
in \cite{Lu1} and were extended to all $p\in \bbR$ in
\cite{WY2008}. Comparing the condition on $K$ in Corollary
\ref{Lp:p:positive:symmetrization} with those in \cite{WY2008}, here 
one does not require the centroid of $K$ to be at the origin. This
was first noticed in \cite{Zhang2007} by Zhang.

 \noindent {\bf Proof.} Let $p\in (0,1)$ and
 $\phi(t)=t^{\frac{p}{n+p}}$ for $t\in (0, \infty)$. Then $\phi\in Conc(0,
\infty)$. Moreover, it is easily checked that
$$F(t)=\phi(t^{n+1})=t^{\frac{np+p}{n+p}}, \ \ t\in (0, \infty)$$ is
strictly concave since $0<\frac{np+p}{n+p}<1$. That is, $\phi
(t)=t^{\frac{p}{n+p}}$ verifies conditions on Theorems
\ref{symmetrization:L:phi} and \ref{L:phi:affine isoperimetric
inequality}. Therefore, (i) follows immediately from Theorem
\ref{symmetrization:L:phi}.

  \noindent For (ii), one first has, by Theorem
\ref{L:phi:affine isoperimetric inequality}, $as_p(K)\leq
as_p(B_K).$ Note that, $B_K=r B_2^n$ with $$r
=\left(\frac{|K|}{|B^n_2|}\right)^{1/n}$$ and for all $\lambda>0$,
$$as_p(\lambda K)=\lambda^{\frac{n(n-p)}{n+p}}as_p(K).$$ Then, one
has $$\frac{as_p(K)}{as_p(B_2^n)}\leq
\frac{as_p(B_K)}{as_p(B_2^n)}=\left(\frac{|K|}{|B^n_2|}
\right)^{\frac{n-p}{n+p}}.$$
As $F$ is strictly concave, by Theorem
\ref{L:phi:affine isoperimetric inequality}, one gets that among all convex bodies with positive Gaussian curvature almost everywhere, the $L_p$ affine surface area attains its maximum {\it only} at ellipsoids.

\subsection{$L_{\psi}$ affine surface areas are decreasing under the
Steiner symmetrization}
For the $L_{\psi}$ affine surface areas, one has the
following theorem.

 \bt\label{symmetrization:L:psi} Let $K\subset
\bbR^n$ be a convex body with the origin in its interior and
$\psi\in Conv(0, \infty)$. Then, if the function
$G(t)=\psi(t^{n+1})$ for $t\in (0, \infty)$ is convex, one has, for all
$\xi \in S^{n-1}$,
\begin{equation*}as_{\psi} (K)\geq as_{\psi}(S_{\xi}(K)).
\end{equation*}  \et

\noindent {\bf Remark.}  Note that the function $\psi \in Conv(0,\infty)$ is monotone
decreasing, hence $G$ is also a decreasing function.   In fact,  in view of Lemma \ref{curvature:symmetrization:1},  the condition that $G$ is convex and monotone decreasing is more natural for proving Theorem \ref{symmetrization:L:psi}. If (non-constant) function $G$ is monotone decreasing and convex with $\lim_{t\rightarrow 0} G(t)=\infty$ and $\lim_{t\rightarrow \infty} G(t)=0$, then $\psi \in Conv(0,\infty)$. (It is easily checked that $G$ constant implies $\psi$ constant).  To this end, it is easy to see that $\lim_{t\rightarrow 0} \psi(t)=\infty$ and $\lim_{t\rightarrow \infty} \psi(t)=0$. For all $0<t<s$, $\psi(t)=G(t^{\frac{1}{n+1}})\geq G(s^{\frac{1}{n+1}})=\psi(s)$ and hence $\psi$ is monotone decreasing. Note that the function $t^{\frac{1}{n+1}}$ is concave, and by $G$ being decreasing,  for all $\lambda\in [0,1]$ and $0<t<s$, 
\begin{eqnarray*} \psi(\lambda t+(1-\lambda) s)&=&G([\lambda t+(1-\lambda) s]^{\frac{1}{n+1}})\leq G(\lambda t^{\frac{1}{n+1}}+(1-\lambda) s^{\frac{1}{n+1}})\\ &\leq & \lambda G(t^{\frac{1}{n+1}})+(1-\lambda) G(s^{\frac{1}{n+1}})=\lambda \psi(t)+(1-\lambda)\psi(s).\end{eqnarray*}  For homogeneous function $\psi(t)=t^a$, to have $G$ convex and monotone decreasing, one needs $a\leq 0$.

 \noindent {\bf Proof of Theorem \ref{symmetrization:L:psi}.} The proof of Theorem
\ref{symmetrization:L:psi} is similar to that of Theorem
\ref{symmetrization:L:phi}. Here, for completeness, we include its
proof with modification emphasized.

Without loss of generality, we assume that $K$ has positive Gaussian curvature almost everywhere. Otherwise, if  $\mu_K(\{y\in \partial K: \kappa_K(y)=0\})>0$, then $as_{\psi}(K)=\infty$, and hence the desired result follows. 

As $as_{\psi}(K)$ is $SL(n)$-invariant, without loss of
generality, we only work on the direction $\xi =e_n=(0, \cdots, 0
,1)$.  Let
$h(x)=[f(x)+g(x)]/2$. Note that $\psi\in Conv(0, \infty)$ is a
decreasing function, so is $G(t)=\psi(t^{n+1})$ on $t\in (0,
\infty)$. By inequality (\ref{determinant:inequality:01}), one has
\begin{eqnarray} \psi\left(\frac{|\mathrm{det}
(d^2h(x))|}{\langle h(x)\rangle^{n+1}}\right)&=&
G\left(\frac{2|\mathrm{det} (d^2h(x))|^{\frac{1}{n+1}}}{\langle
f(x)\rangle+\langle g(x)\rangle}\right)\nonumber\\ &\leq&
G\left(\!\frac{\left|\mathrm{det} (d^2
f(x))\right|^{\frac{1}{n+1}} \!+\!\left|\mathrm{det} (d^2
g(x))\right|^{\frac{1}{n+1}}}{\langle f(x)\rangle+\langle
g(x)\rangle}\!\right)\nonumber\\
&\leq& G\left(\frac{|\mathrm{det} (d^2
f(x))|^{\frac{1}{n+1}}}{\langle f(x)\rangle} \right)\frac{\langle
f(x)\rangle}{\langle f(x)\rangle+\langle g(x)\rangle}\nonumber\\
&& \ + \ \ G\left(\!\frac{\left|\mathrm{det} (d^2
g(x))\right|^{\frac{1}{n+1}}}{\langle
g(x)\rangle}\!\!\right)\frac{\langle g(x)\rangle}{\langle
f(x)\rangle+\langle g(x)\rangle}, \label{convexity:inequality:1}
\end{eqnarray} where inequality (\ref{convexity:inequality:1}) follows from the
convexity of the function $G(t)$ on $(0, \infty)$. Therefore, by
(\ref{symmetrization:en}) and Lemma
\ref{Equivalent:definition:affine:surface:area:psi}, we have for
all $\psi\in Conv(0, \infty)$,
\begin{eqnarray*}
as_{\psi}(S_{e_n}(K))&=&2 \int
_{K_0}\left\{\psi\left(\frac{|\mathrm{det} (d^2h(x))|}{\langle
h(x)\rangle^{n+1}}\right)
\langle h(x)\rangle\right\}\,dx \nonumber\\
&\leq& \int _{K_0}\left\{\psi\left(\frac{|\mathrm{det}
(d^2f(x))|}{\langle f(x)\rangle^{n+1}}\right) \langle
f(x)\rangle+\psi\left(\frac{|\mathrm{det} (d^2g(x))|}{\langle
g(x)\rangle^{n+1}}\right) \langle g(x)\rangle\right\}\,dx
\nonumber
\\&=&as_{\psi}(K). \end{eqnarray*}

 Let $K$ be a convex body with curvature function. Similar 
 to the $L_{\phi}^*$ affine surface
area,  the $L_{\psi}^*$ affine surface area for $\psi\in
Conv(0, \infty)$ can be formulated as 
$$as_{\psi}^*(K)=\int_{S^{n-1}}\psi(f_{-n}(K,u))\,d\nu_K(u).$$
Notice that the $L_p$ affine surface area for $p<-n$ is a special
case for the $L_{\psi}^*$ affine surface area with
$\psi(t)=t^{\frac{n}{n+p}}$. It was prove that
\begin{equation} \label{duality:convex}
as_{\psi}^*(K)=as_{\psi}(K^\circ),\end{equation} for all convex
body $K$ having curvature function and with the origin in its interior \cite{Ludwig2009}.
Combining with Theorem \ref{symmetrization:L:psi}, one immediately
has the following result.

\bc Let $K\subset \bbR^n$ be a convex body having curvature function and with the origin in its
interior. Let $\psi\in Conv(0, \infty)$. Assume that the function
$G(t)=\psi(t^{n+1})$ for $t\in (0, \infty)$ is convex. Then, the
$L_{\psi}^*$ affine surface area is monotone decreasing under the
Steiner symmetrization in the following sense:
$$as_{\psi}^* (K)\geq as_{\psi}^*([S_{\xi}(K^\circ)]^\circ),$$ for
all $\xi\in S^{n-1}$ such that,  $[S_{\xi}(K^\circ)]^\circ$ has curvature function. \ec

\bt \label{L:psi:affine isoperimetric inequality} Let $K$ be a
convex body with the origin in its interior, and $B_K$ be the
origin-symmetric ball such that $|K|=|B_K|$. Then, for all $\psi
\in Conv(0, \infty)$ such that the function $G(t)=\psi(t^{n+1})$
for $t\in (0, \infty)$ is convex, one has
\begin{equation*}  as_{\psi}(K)\geq
as_{\psi}(B_{K}).\end{equation*} If in addition $G(t)$ is strictly convex, equality holds if and only if $K$ is an origin-symmetric 
ellipsoid. \et

 \noindent {\bf Proof.} The proof of Theorem
\ref{L:psi:affine isoperimetric inequality} is almost identical to
that of Theorem \ref{L:phi:affine isoperimetric inequality}. Here,
we only mention the main modification.

Let $K\subset \bbR^n$ be a convex body with the origin
in its interior. Suppose that $\psi\in Conv(0, \infty)$. As in the
proof of Theorem \ref{L:phi:affine isoperimetric inequality}, one can find  a sequence of directions   $\{u_i\}_{i=1}^{\infty}\subset \Omega$ such that  $K_i$ converges to
$B_K$ in the Hausdorff distance.  Here $K_i$ is defined as follows:  
$$K_1=S_{u_1}(K); \ \ \  K_{i+1}=S_{u_{i+1}}(K_i), \forall i=1, 2, \cdots $$ 
  Theorem \ref{symmetrization:L:psi} implies that
$$as_{\psi}(K)\geq as_{\psi}(K_1)\geq \cdots \geq as_{\psi}(K_j), \ \ \forall
j\in \bbN. $$ Combining with the lower--semicontinuity of
$as_{\psi}(\cdot)$, one has,
\begin{eqnarray*} as_{\psi}(K)&\geq& \lim
\inf _{j\rightarrow \infty} as_{\psi}(K_j)\geq as_{\psi}(\lim
_{j\rightarrow \infty} K_j)=as_{\psi}(B_K). \end{eqnarray*}

Now let us assume that $G$ is strictly convex and $as_{\psi}(K)=as_{\psi}(B_K)$. Clearly, to have $as_{\psi}(K)=as_{\psi}(B_K)$, $K$ must have positive Gaussian curvature a.e. on $\partial K$, as otherwise $as_{\psi}(K)=\infty$.

 Let $K$ be a convex body with positive Gaussian curvature almost everywhere. We now claim that the set $M(K, e_n)$ is contained in a hyperplane. In this case, we assume that $e_n$ is a direction such that  both the ovegraph and undergraph functions $f, g$ are differentiable  at $0$.  Equation $as_{\psi}(K)=as_{\psi}(S_{e_n}(K))$ requires equalities for (\ref{convexity:inequality:1}). By the strict convexity of $G$, one has, for almost every $x\in  K_0$, 
 \begin{equation}\label{representation:affine:invariant-1} \left|\mathrm{det}
(d^2 f(x))\right|=\left|\mathrm{det} (d^2
g(x))\right|;  \ \ \  \frac{|\mathrm{det} (d^2
f(x))|^{\frac{1}{n+1}}}{\langle f(x)\rangle} =\frac{\left|\mathrm{det} (d^2
g(x))\right|^{\frac{1}{n+1}}}{\langle
g(x)\rangle}>0.\end{equation} Hence, $f(x)-g(x)=\langle x, \nabla (f(x)-g(x)) \rangle$ for almost all $x\in K_0$.  Assume that both $f, g$ are differentiable at $0$ . From Lemma 4.3 in \cite{Hug1995}, one obtains that $f(x)-g(x)$ is linear, and hence $M(K, e_n)$ is contained in a hyperplane. 

Let $\Omega$ be the dense subset of $S^{n-1}$ such that the corresponding  overgraph and undergraph functions are both differentiable at $0$.  For any $u\in \Omega$, there is a rotation $T$ such that $T(u)$ is parallel to $e_n$. The above claim then implies that $M(TK, T(u))$ (and hence $M(K, u)$) is contained in a hyperplane. By Lemma  \ref{Charactristic:Ellipsoid:Hug}, $K$ is an ellipsoid.  Moreover, $K$ has to be an origin-symmetric ellipsoid. To this end,  we assume that $K$ is a ball with center $y_c\neq 0$. By formulas (\ref{representation:affine:invariant}) and (\ref{representation:affine:invariant-1}), one gets 
 $$\frac{[\kappa_K(y)]^{\frac{1}{n+1}}} {\langle y,N_{K}(y)\rangle
}=\frac{|\mathrm{det} (d^2
f(x))|^{\frac{1}{n+1}}}{\langle f(x)\rangle} =\frac{\left|\mathrm{det} (d^2
g(x))\right|^{\frac{1}{n+1}}}{\langle
g(x)\rangle}=\frac{[\kappa_K(z)]^{\frac{1}{n+1}}} {\langle z,N_{K}(z)\rangle}, \ \ \forall x\in K_{e_n},$$ with $y=(x, f(x))\in \partial K$ and $z=(x, -g(x))\in \partial K$. Note that the curvature of $K$ is a constant, and hence $\kappa_K(y)=\kappa_K(z)$. This implies  $\langle y,N_{K}(y)\rangle=\langle z,N_{K}(z)\rangle$, a contradiction with $K$ being a ball with center $y_c\neq 0$.

Combining with formula (\ref{duality:convex}) one
immediately has the following isoperimetric inequality for the
$L_{\psi}^*$ affine surface area.

\bc Let $K\subset \bbR^n$ be a convex body having curvature function and with the origin in its
interior. Let $B_K$ be the origin-symmetric Euclidean ball with
$|K|=|B_K|$. For all $\psi \in Conv(0, \infty)$ such that the
function $G(t)=\psi(t^{n+1})$ for $t\in (0, \infty)$ is convex,
one has
\begin{equation*}as_{\psi}^*(K)\geq
as_{\psi}^*([B_{K^\circ}]^\circ).\end{equation*} If in addition  $G(\cdot)$ is strictly convex, equality holds if
and only if $K$ is an origin-symmetric ellipsoid.\ec

 The $L_p$ affine surface areas for
$p\in (-n, 0)$ are special cases of the $L_{\psi}$ affine surface
areas with $\psi(t)=t^{\frac{p}{n+p}}$. We have the following
results.
 \bc \label{Lp:p:negative:symmetrization} Let $K$ be a convex body with
 the origin in its interior, and
 let $p\in (-n, 0)$.\\
 \noindent (i) The $L_p$ affine surface area for $p\in (-n, 0)$
 is monotone decreasing under the Steiner symmetrization.
 That is,  for any $\xi \in S^{n-1}$, one has
\begin{eqnarray*}as_p(K)\geq as_p(S_{\xi}(K)).
 \end{eqnarray*} (ii) The $L_p$
affine surface areas attain their minimum at the ellipsoid, among all convex bodies with fixed volume. More
precisely, \begin{eqnarray*} \frac{as_p(K)}{as_p(B^n_2)}\geq
\left(\frac{|K|}{|\bn n|}\right)^{\frac{n-p}{n+p}}.
\end{eqnarray*} Equality holds if and only if $K$ is an origin-symmetric ellipsoid.
 \ec

 \noindent {\bf Remark.} The $L_p$ affine isoperimetric
inequalities for $p\in (-n,0)$ were first established in
\cite{WY2008}. Comparing the condition on $K$ in Corollary
\ref{Lp:p:negative:symmetrization} with those in \cite{WY2008}, here 
one does not require the centroid of $K$ to be at the origin. This
was first noticed in \cite{Zhang2007}.

 \noindent {\bf Proof.} Let $p\in (-n, 0)$ and
 $\psi(t)=t^{\frac{p}{n+p}}$ for $t\in (0, \infty)$. Then $\psi\in Conv(0,
\infty)$. Moreover, it is easily checked that
$$G(t)=\psi(t^{n+1})=t^{\frac{np+p}{n+p}}, \ \ t\in (0, \infty)$$ is convex since
$\frac{np+p}{n+p}<0$. That is, $\psi (t)=t^{\frac{p}{n+p}}$
verifies conditions on Theorems \ref{symmetrization:L:psi} and
\ref{L:psi:affine isoperimetric inequality}. Therefore, (i)
follows immediately from Theorem \ref{symmetrization:L:psi}.

 For (ii), one first has, by Theorem \ref{L:psi:affine
isoperimetric inequality}, $as_p(K)\geq as_p(B_K).$ Note that,
$B_K=r B_2^n$ with $$r =\left(\frac{|K|}{|B^n_2|}\right)^{1/n}$$
and for all $\lambda>0$,
$$as_p(\lambda K)=\lambda^{\frac{n(n-p)}{n+p}}as_p(K).$$ Then, one
has $$\frac{as_p(K)}{as_p(B_2^n)}\geq
\frac{as_p(B_K)}{as_p(B_2^n)}=\left(\frac{|K|}{|B^n_2|}
\right)^{\frac{n-p}{n+p}}.$$
As $G$ is strictly convex, by Theorem \ref{L:psi:affine
isoperimetric inequality}, one gets that equality holds if and only if $K$ is an origin-symmetric ellipsoid.

 \noindent {\bf Acknowledgments.} The author is grateful to the referee. This paper is supported
 by a NSERC grant and a start-up grant from the Memorial University
 of Newfoundland.

\vskip 2mm \noindent Deping Ye, \ \ \ {\small \tt deping.ye@mun.ca}\\
{\small \em Department of Mathematics and Statistics\\
   Memorial University of Newfoundland\\
   St. John's, Newfoundland, Canada A1C 5S7 }
\end{document}